\documentclass[12pt, leqno]{amsart}

\usepackage{graphicx}
\usepackage{amssymb,amsmath,amsthm,amscd}
\usepackage{mathrsfs}
\usepackage{enumerate}
\usepackage[usenames,dvipsnames]{color}
\usepackage[colorlinks=true, pdfstartview=FitV,
 linkcolor=blue,citecolor=blue,urlcolor=blue]{hyperref}

\usepackage[all]{xy}

\allowdisplaybreaks[3]

\newcommand{\nc}{\newcommand}
\numberwithin{equation}{section}

\theoremstyle{plain}
\newtheorem{lemma}{Lemma}[section]
\newtheorem{prop}[lemma]{Proposition}
\newtheorem{theorem}[lemma]{Theorem}
\newcommand{\Prop}{\begin{prop}}
\newcommand{\enprop}{\end{prop}}
\newcommand{\Lemma}{\begin{lemma}}
\newcommand{\enlemma}{\end{lemma}}
\newcommand{\Th}{\begin{theorem}}
\newcommand{\enth}{\end{theorem}}
\newtheorem{corollary}[lemma]{Corollary}
\newcommand{\Cor}{\begin{corollary}}
\newcommand{\encor}{\end{corollary}}
\newtheorem{definition}[lemma]{Definition}
\newtheorem*{conjecture}{Conjecture}
\newcommand{\Def}{\begin{definition}}
\newcommand{\edf}{\end{definition}}
\newtheorem{sublemma}[lemma]{Sublemma}
\newcommand{\Sublemma}{\begin{sublemma}}
\newcommand{\ensub}{\end{sublemma}}

\theoremstyle{definition}
\newtheorem{remark}[lemma]{Remark}

\newtheorem{Convention}[lemma]{Convention}
\newcommand{\Conv}{\begin{Convention}}
\newcommand{\enconv}{\end{Convention}}
\nc{\Con}{\begin{conjecture}}
\nc{\encon}{\end{conjecture}}
\nc{\Rem}{\begin{remark}}
\nc{\enrem}{\end{remark}}
\newcommand{\C}{{\mathbb C}}
\newcommand{\Q}{\mathbb {Q}}

\newcommand{\Z}{{\mathbb Z}}
\newcommand{\B}{{\mathbf{B}}}

\newcommand{\one}{{\bf{1}}}
\newcommand{\seteq}{\mathbin{:=}}

\newcommand{\g}{{\mathfrak{g}}}

\newcommand{\Hom}[1][\mspace{-.5mu}]{\mathrm{Hom}_{\raise1.5ex\hbox to.1em{}%
\mspace{.5mu} #1}}
\newcommand{\End}{\operatorname{End}}

\newcommand{\isoto}[1][]{\mathop{\xrightarrow%
[{\raisebox{.3ex}[0ex][.3ex]{$\scriptstyle{#1}$}}]%
{{\raisebox{-.6ex}[0ex][-.6ex]{$\mspace{2mu}\sim\mspace{2mu}$}}}}}

\newcommand{\eq}{\begin{eqnarray}}
\newcommand{\eneq}{\end{eqnarray}}

\newcommand{\hs}{\hspace*}

\newcommand{\To}[1][{\hs{2ex}}]{\xrightarrow{\,#1\,}}

\newcommand{\eqn}{\begin{eqnarray*}}
\newcommand{\eneqn}{\end{eqnarray*}}
\newcommand{\on}{\operatorname}
\newcommand{\Ker}{\on{Ker}}

\newcommand{\bna}{\be[{\rm(a)}]}

\newcommand{\QED}{\end{proof}}
\newcommand{\Proof}{\begin{proof}}

\newcommand{\soplus}{\mathop{\mbox{\normalsize$\bigoplus$}}\limits}

\newcommand{\id}{\on{id}}
\newcommand{\ba}{\begin{array}}
\newcommand{\ea}{\end{array}}

\newcommand{\monoto}{\rightarrowtail}

\newcommand{\set}[2]{\left\{#1 \mid #2 \right\}}
\newcommand{\supp}{\operatorname{supp}}

\newcommand{\eqsub}{\begin{subequations}\begin{eqnarray}}
\newcommand{\eneqsub}{\end{eqnarray}\end{subequations}}

\newcommand{\ol}{\overline}

\nc{\la}{\lambda}
\nc{\lam}{\lambda}
\nc{\U}[1][\g]{U_q(#1)}
\nc{\te}{\tilde{e}}
\nc{\tei}{\tilde{e}_i}
\nc{\tf}{\tilde{f}}
\nc{\tfi}{\tilde{f}_i}
\nc{\tU}{\widetilde U_q(\g)}
\nc{\tE}{\tilde{E}}
\nc{\tF}{\widetilde{\F}}
\nc{\tK}{\widetilde{K}}

\nc{\tk}{\tilde{k}}
\nc{\tkone}{\tk_{\ol{1}}}
\nc{\teone}{\tilde{e}_{\ol{1}}}
\nc{\tfone}{\tilde{f}_{\ol{1}}}

\nc{\teibar}{\tilde{e}_{\ol{i}}} \nc{\tfibar}{\tilde{f}_{\ol{i}}}
\nc{\tki}{{\tk}_{\ol {i}}}

\nc{\BZ}{{\mathbb{Z}}}
\nc{\al}{\alpha}
\nc{\qs}{{q}}
\nc{\lan}{\langle}
\nc{\ran}{\rangle}
\nc{\re}{{\mathrm{re}}}
\nc{\wt}{\operatorname{wt}}
\nc{\ch}{\operatorname{ch}}
\nc{\Uf}[1][\g]{U^-_q(#1)}
\nc{\Ue}{U^+_q(\g)}
\nc{\eps}{\varepsilon}
\nc{\vphi}{\varphi}
\nc{\sphi}{\varphi^*}
\nc{\seps}{\varepsilon^*}
\newcommand{\bB}{\mathbf{B}}
\nc{\nn}{\nonumber}

\nc{\vp}{\varpi}
\nc{\cls}{{\operatorname{cl}}}
\nc{\Wt}{{\operatorname{Wt}}}
\nc{\Us}{U'_q(\g)}
\nc{\La}{\Lambda}
\nc{\ro}{{\rm(}}
\nc{\rf}{{\rm)}}
\nc{\norm}{{\mathrm{norm}}}
\nc{\qbox}{\quad\mbox}
\nc{\braid}{{\mathfrak{B}}}
\nc{\Ad}{\operatorname{Ad}}
\nc{\Aut}{\operatorname{Aut}}
\nc{\dt}[1]{\tilde{\tilde #1}}
\nc{\Sn}{S^{{\mathrm{norm}}}}
\nc{\aff}{{\rm{aff}}}
\nc{\rk}{{\mathrm{rk}}}
\nc{\tP}{\widetilde{P}}
\nc{\tW}{\widetilde{W}}
\nc{\Dyn}{\mathrm{Dyn}}
\nc{\tD}{\widetilde{\Delta}}
\nc{\height}[1]{{\operatorname{ht}}(#1)}
\nc{\bl}{\bigl(}
\nc{\br}{\bigr)}
\nc{\Hecke}{\mathrm{H}}
\nc{\HA}{\Hecke^{\mathrm{A}}}
\nc{\HB}{\Hecke^{\mathrm{B}}}
\newcommand{\scbul}{{\,\raise1pt\hbox{$\scriptscriptstyle\bullet$}\,}}
\nc{\vac}{{\phi}}
\nc{\Bt}{\B_\theta(\g)}
\nc{\be}{\begin{enumerate}}
\nc{\ee}{\end{enumerate}}
\nc{\low}{{\mathrm{low}}}
\nc{\upper}{{\mathrm{up}}}
\nc{\Zodd}{\Z_{\mathrm{odd}}}
\nc{\Ft}[1][n]{\mathbb{P}\mathrm{ol}_{#1}}
\nc{\Ftf}[1][n]{\widetilde{\mathbb{P}\mathrm{ol}}_{#1}}
\nc{\KA}{\on{K}^{\mathrm{A}}}
\nc{\KB}{\on{K}^{\mathrm{B}}}
\nc{\Res}{\on{Res}}
\nc{\Fc}[1][{n,m}]{\mathbf{F}_{#1}}
\nc{\tphi}{\tilde{\varphi}}
\nc{\CO}{\mathscr{O}}
\nc{\inte}{\mathrm{int}}
\nc{\Oint}{\mathcal{O}^{\ge0}_{\inte}}
\nc{\vs}{\vspace*}
\nc{\tL}{\widetilde{L}}
\nc{\tu}{\tilde{u}}
\nc{\noi}{\noindent}
\nc{\heigh}{\mathfrak{t}}
\nc{\lowest}{\mathfrak{l}}
\nc{\rootl}{\mathsf{Q}}
\nc{\cl}{\mathrm{cl}}
\nc{\uqpg}{U'_q(\mathfrak g)}
\nc{\uq}{\uqpg}
\nc{\Oh}{\widehat{\mathcal{O}}}

\nc{\KLR}{quiver Hecke algebra}
\nc{\KLRs}{quiver Hecke algebras}
\nc{\cor}{\mathbf{k}}
\nc{\cora}{{\cor(A)}}
\nc{\haut}{\mathrm{ht}}
\nc{\tens}{\mathop\otimes}
\nc{\gmod}{\mbox{-$\mathrm{mod}$}}
\nc{\proj}{\mbox{-$\mathrm{proj}$}}
\nc{\gproj}{\mbox{-$\mathrm{gproj}$}}
\nc{\smod}{\mbox{-$\mathrm{mod}$}}

\nc{\h}{\mathfrak h}
\nc{\Rnorm}{R^{\rm{norm}}}
\nc{\K}{\C(q)}
\nc{\Vhat}{\widehat{V}}
\nc{\F}{\mathcal{F}}

\def\T{{\mathcal T}}

\nc{\fd}[1][A]{\on{\mathrm{flat.dim}_{#1}}}
\nc{\bP}{{\mathbb{P}}}
\nc{\bPh}{\widehat{\mathbb{P}}}
\nc{\bK}[1][{n}]{\widehat{\mathbb{K}}_{#1}}
\nc{\bV}[1][{n}]{\widehat{V}^{\otimes{#1}}}
\nc{\bVK}[1][{n}]{\widehat{V}^{\otimes{#1}}_{\widehat{\mathbb{K}}}}
\nc{\hV}{\widehat{V}}
\nc{\opp}{\mathrm{opp}}
\nc{\col}{\colon}
\nc{\bnum}{\be[{\rm(i)}]}
\nc{\oep}{\epsilon}
\nc{\qtext}{\quad\text}
\nc{\qtextq}[1]{\quad\text{#1}\quad}
\nc{\longtwoheadrightarrow}[1][]{\xymatrix{\ar@{->>}[r]^-{{#1}}&}}
\nc{\epiTo}[1][]{\longtwoheadrightarrow[{#1}]}
\nc{\epito}{\twoheadrightarrow}
\nc{\monoTo}[1][]{\xymatrix{\ar@{>->}[r]^-{{#1}}&}}
\nc{\sym}{\mathfrak{S}}
\nc{\inp}[1]{{({#1})_{\mathrm{n}}}}
\nc{\rtl}{\rootl}
\nc{\wtd}{\widetilde}
\nc{\etens}{\boxtimes}
\nc{\ds}[1]{\mathrm{d}(#1)}
\nc{\rmat}[1]{{\mathbf{r}}_%
{\mspace{-2mu}\raisebox{-.6ex}{${\scriptstyle{#1}}$}}}
\nc{\rmats}[1]{{\mathbf{r}}_%
{\mspace{-2mu}\raisebox{-.6ex}{${\scriptscriptstyle{#1}}$}}}
\nc{\shc}{\mathcal{C}}
\nc{\Fct}{{\on{Fct}}}
\nc{\tC}{\widetilde{\shc}}
\nc{\Zp}{\Z_{\ge0}}
\nc{\tPhi}{\widetilde{\Phi}}
\nc{\tT}{{\widetilde{\T}}}
\nc{\Ob}{\on{Ob}}
\nc{\bwr}{\mbox{\large$\wr$}}
\nc{\Img}{\on{Im}}
\nc{\Ab}{\mathcal{A}^{\mathrm{big}}}
\nc{\Sb}{\mathcal{S}^{\mathrm{big}}}
\nc{\As}{\mathcal{A}}
\nc{\Ss}{\mathcal{S}}
\nc{\ntens}{\widetilde{\otimes}}
\nc{\hR}{\widehat{R}}
\nc{\nconv}{\star}
\nc{\ts}{\tilde{s}}
\nc{\sho}{\mathcal{O}}
\nc{\bc}{\begin{cases}}
\nc{\ec}{\end{cases}}
\nc{\slnh}{{\widehat{\mathfrak{sl}}_N}}
\nc{\UA}{U_q'(\slnh)}
\nc{\KR}{R_K}
\nc{\cQ}{\mathcal{Q}}
\nc{\Irr}{\mathcal{I}rr}
\nc{\tQ}{\widetilde{\cQ}}
\nc{\bs}{\mathbf{s}}
\nc{\bL}{\mathbb{L}}
\nc{\tg}{\tilde{g}}
\nc{\conv}{\mathbin{\mbox{\large $\circ$}}}
\nc{\hconv}{\mathbin{\mbox{\large $\diamond$}}}

\nc{\bQ}{\ol{Q}}
\renewcommand{\Im}{\on{Im}}
\nc{\xmono}{\ar@{>->}}
\nc{\xepi}{\ar@{->>}}
\nc{\db}[1]{\raisebox{-.5ex}[2ex][1.8ex]{$#1$}}
\nc{\wb}[1]{\mbox{$\rule[-1.1ex]{0ex}{2ex}#1$}}
\nc{\univ}{\mathrm{univ}}
\nc{\rM}{{}^*\mspace{-2mu}M}
\nc{\lM}{M^*}
\nc{\uqm}{\uq\smod}
\nc{\tR}{\widetilde{R}_{\gamma,\beta}}
\nc{\tx}{\tilde{x}}
\nc{\bi}{\mathbf{i}}
\nc{\bj}{\mathbf{j}}
\nc{\ttau}{\widetilde{\tau}}

\newlength{\mylength}
\title
{Simplicity of heads and socles of tensor products}

\author[S.-J. Kang, M. Kashiwara, M. Kim, Se-jin Oh]{Seok-Jin Kang$^{1}$, Masaki Kashiwara$^{2}$, Myungho Kim, Se-jin Oh$^3$}

\address{Department of Mathematical Sciences
         and
         Research Institute of Mathematics \\
         Seoul National University \\ Seoul 151-747, Korea}

         \email{sjkang@snu.ac.kr}

\address{Research Institute for Mathematical Sciences \\
          Kyoto University \\ Kyoto 606-8502, Japan \\
          \& Department of Mathematical Sciences
         and
         Research Institute of Mathematics \\
         Seoul National University \\ Seoul 151-747, Korea}

         \email{masaki@kurims.kyoto-u.ac.jp}

\address{School of Mathematics, Korea Institute for Advanced Study \\ Seoul 130-722, Korea}
         \email{mhkim@kias.re.kr}

\address{Department of Mathematical Sciences
         and
         Research Institute of Mathematics \\
         Seoul National University \\ Seoul 151-747, Korea}
         \email{sj092@snu.ac.kr}

\thanks{$^1$ This work was supported by NRF grant \# 2014021261
and NRF grant \# 2013055408
}
\thanks{$^2$ This work was supported by Grant-in-Aid for
Scientific Research (B) 22340005, Japan Society for the Promotion of
Science.}
\thanks{$^3$ This work was supported by BK21 PLUS SNU Mathematical Sciences Division}

\keywords{Quantum affine algebra, Khovanov-Lauda-Rouquier algebra,
$R$-matrix}

\subjclass[2010]
{81R50, 16G, 17B37}
\date{April 15, 2014}

\begin{document}

\begin{abstract}
We prove that,  for simple modules $M$ and $N$ over  a quantum affine algebra,
their tensor product $M\tens N$ has a simple head and a simple socle
if $M\tens M$ is simple.
A similar result is proved for the convolution product of
simple modules over \KLRs.

\end{abstract}

\maketitle

\section*{Introduction}

Let $\mathsf{g}$ be a complex simple Lie algebra and $U_q(\mathsf{g})$  the associated
quantum group. The multiplicative property of the upper global
basis $\bB$ of the negative half $U_q^-(\mathsf{g})$ was investigated in
\cite{BZ93,L03}. Set $q^\Z\bB=\set{q^nb}{b\in\bB,\,n\in\Z}$. In
\cite{BZ93}, Berenstein and Zelevinsky conjectured that, for $b_1$,
$b_2\in\bB$, the product $b_1b_2$ belongs to $q^\Z\bB$ if and only
if $b_1$ and $b_2$ $q$-commute, (i,e.,\ $b_2b_1=q^nb_1b_2$ for some
$n\in\Z$). However, Leclerc found examples of $b\in\bB$ such that
$b^2\not\in q^\Z\bB$ (\cite{L03}).

On the other hand, the algebra $U_q^-(\mathsf{g})$ is categorified by \KLRs\ 
(\cite{KL09, KL11, R08}) and also by quantum affine
algebras (\cite{HL10,HL11,K^3,K^3II}). In this context,
the products in  $U_q^-(\mathsf{g})$ correspond to the convolution
or the tensor products in \KLRs\ or quantum affine algebras. 
The upper global basis corresponds to the set of isomorphism classes of simple
modules over the \KLR\ or the quantum affine
algebras (\cite{Ariki,R11,VV09}) under suitable conditions.
Then Leclerc conjectured
several properties of products of upper global bases and also
 convolutions and  tensor products of simple modules. The
purpose of this paper is to give an affirmative answer to some of
his conjectures.

In this introduction, we state our results
in the case of modules over quantum affine algebras.
The similar results hold also for \KLRs\ (see \S\,\ref{sec:qHal}).

\medskip
Let $\g$ be an affine Lie algebra and $\uqpg$ 
the associated quantum affine algebra. A simple
$\uq$-module $M$ is called {\em real} if $M\tens M$ is also simple.

\Con[{\cite[Conjecture 3]{L03}}]
Let $M$ and $N$ be finite-dimensional simple $\uqpg$-modules.
We assume further
that $M$ is real. Then
$M\tens N$ has a simple socle $S$ and a simple head $H$.
Moreover, if $S$ and $H$ are isomorphic, then $M\tens N$ is simple.
\encon

In this paper, we shall give an affirmative answer to
this conjecture (Theorem~\ref{th:main A} and Corollary~\ref{cor:mainQ}). In the course of
the proof, $R$-matrices play an important role.
Indeed, the simple socle of $M\tens N$ coincides with the
image of the renormalized $R$-matrix
$\rmat{N,M}\col N\tens M\to M\tens N$
and the simple head of $M\tens N$ coincides with the
image of the renormalized $R$-matrix
$\rmat{M,N}\col M\tens N\to N\tens M$.

Denoting by $M\hconv N$ the head of $M\tens N$,
we also prove that
$N\mapsto M\hconv N$ is an automorphism of
the set of the isomorphism classes of simple $\uq$-modules
(Corollary~\ref{cor:bj}). The inverse is given by $N\mapsto N\hconv \rM$,
where $\rM$ is the right dual of $M$. It is an analogue of 
Conjecture 2 in \cite{L03} originally stated for global bases.

\bigskip
\noi
{\bf Acknowledgements.} We would like to thank
Bernard Leclerc for many fruitful discussions and
his kind explanations on his works.

\section{Quiver Hecke algebras}\label{sec:KLR}

In this section, we briefly recall the basic facts on \KLRs\ and $R$-matrices
following \cite{K^3}. 
Since the grading of \KLRs\ is not important in this paper,
we ignore the grading.
Throughout the paper, {\em modules mean left modules.}

\subsection{Convolutions}
We shall recall the definition of \KLRs.
Let $\cor$ be a field.  Let $I$ be an index set. Let $\rtl$
be the free $\Z$-module with a basis $\{\al_i\}_{i\in I}$. Set
$\rtl^+=\sum_{i\in I}\Z_{\ge0}\al_i$. For
$\beta=\sum_{k=1}^n\al_{i_k}\in\rtl^+$, we set $\height{\beta}=n$.
For $n \in \Z_{\ge 0}$ and $\beta \in \rootl^+$ such that
$\height{\beta} = n$, we set
$$I^{\beta} = \set{\nu = (\nu_1, \ldots, \nu_n) \in I^{n}}%
{ \alpha_{\nu_1} + \cdots + \alpha_{\nu_n} = \beta }.$$

Let us take  a family of polynomials $(Q_{ij})_{i,j\in I}$ in $\cor[u,v]$
which satisfy
\eqn
&&Q_{ij}(u,v)=Q_{ji}(v,u)\quad \text{for any $i$, $j\in I$,}\\
&&Q_{ii}(u,v)=0\quad \text{for any $i\in I$.}
\eneqn
For $i,j\in I$, we set
$$\bQ_{ij}(u,v,w)=\dfrac{Q_{ij}(u,v)-Q_{ij}(w,v)}{u-w}\in\cor[u,v,w].$$

We denote by
$\sym_{n} = \langle s_1, \ldots, s_{n-1} \rangle$ the symmetric group
on $n$ letters, where $s_i\seteq (i, i+1)$ is the transposition of $i$ and $i+1$.
Then $\sym_n$ acts on $I^n$ by place permutations.

\begin{definition}
For $\beta \in \rootl^+$ with $\height{\beta}=n$, the {\em
\KLR}  $R(\beta)$  at $\beta$ associated
with a matrix $(Q_{ij})_{i,j \in I}$ is the $\cor$-algebra
generated by the elements $\{ e(\nu) \}_{\nu \in  I^{\beta}}$, $
\{x_k \}_{1 \le k \le n}$, $\{ \tau_k \}_{1 \le k \le n-1}$
satisfying the following defining relations{\rm:}
\begin{equation*} \label{eq:KLR}
\begin{aligned}
& e(\nu) e(\nu') = \delta_{\nu, \nu'} e(\nu), \ \
\sum_{\nu \in  I^{\beta} } e(\nu) = 1, \\
& x_{k} x_{m} = x_{m} x_{k}, \ \ x_{k} e(\nu) = e(\nu) x_{k}, \\
& \tau_{m} e(\nu) = e(s_{m}(\nu)) \tau_{m}, \ \ \tau_{k} \tau_{m} =
\tau_{m} \tau_{k} \ \ \text{if} \ |k-m|>1, \\
& \tau_{k}^2 e(\nu) = Q_{\nu_{k}, \nu_{k+1}} (x_{k}, x_{k+1})
e(\nu), \\
& (\tau_{k} x_{m} - x_{s_k(m)} \tau_{k}) e(\nu) = \begin{cases}
-e(\nu) \ \ & \text{if $m=k$ and $\nu_{k} = \nu_{k+1}$,} \\
e(\nu) \ \ & \text{if $m=k+1$ and $\nu_{k}=\nu_{k+1}$,} \\
0 \ \ & \text{otherwise},
\end{cases} \\
& (\tau_{k+1} \tau_{k} \tau_{k+1}-\tau_{k} \tau_{k+1} \tau_{k}) e(\nu)
=\begin{cases} \bQ_{\nu_{k}, \nu_{k+1}}(x_{k},x_{k+1},x_{k+2})& \text{if
$\nu_{k} = \nu_{k+2}$,} \\
0 \ \ & \text{otherwise}.
\end{cases}
\end{aligned}
\end{equation*}
\end{definition}

For an element $w$ of the symmetric group $\sym_n$,
let us choose  a  reduced expression $w=s_{i_1}\cdots s_{i_\ell}$, and
set $$\tau_w=\tau_{i_1}\cdots \tau_{i_\ell}.$$
In general, it depends on the choice of
reduced expressions of $w$.
Then we have {the PBW decomposition}
\eq
&&R(\beta)=
\soplus_{\nu\in  I^{\beta} , \;w\in\sym_n}\cor[x_1,\ldots, x_n]e(\nu)\tau_w.
\label{eq:PBW}
\eneq
We denote by $R(\beta)\smod$ 
the category of $R(\beta)$-modules $M$ such that
$M$ is finite-dimensional over $\cor$
and the action of $x_k$ on $M$ is nilpotent for any $k$.

\smallskip
For an $R( \beta)$-module $M$, the dual space
$$M^* \seteq \Hom[\cor](M, \cor)$$
is endowed with the $R(\beta)$-module structure given by 
\begin{align*}
(r \cdot  f)(u) \seteq f(\psi(r) u) \quad 
\text{for $f\in M^*$, $r \in R(\beta)$ and $u \in M$,}
\end{align*}
where $\psi$ denotes the $\cor$-algebra anti-involution on $R(\beta
)$ which fixes the generators $\{ e(\nu) \}_{\nu \in  I^{\beta}}$, $
\{x_k \}_{1 \le k \le n}$, $\{ \tau_k \}_{1 \le k \le n-1}$.

\medskip
For $\beta, \gamma \in \rootl^+$ with $\height{\beta}=m$ and $\height{\gamma}= n$,
 set
$$e(\beta,\gamma)=\displaystyle\sum_{\substack%
{\nu \in I^{m+n}, \\ (\nu_1, \ldots ,\nu_m) \in I^{\beta},\\
(\nu_{m+1}, \ldots ,\nu_{m+n}) \in I^{\gamma}}} e(\nu) \in R(\beta+\gamma). $$
Then $e(\beta,\gamma)$ is an idempotent.
Let
\eqn R( \beta)\tens R( \gamma  )\to e(\beta,\gamma)R( \beta+\gamma)e(\beta,\gamma) \label{eq:embedding}
\eneqn
be the $\cor$-algebra homomorphism given by 
\eqn
&e(\mu)\tens e(\nu)\mapsto e(\mu*\nu)  & (\mu\in I^{\beta}, \
 \nu\in I^{\gamma}) \\
&x_k\tens 1\mapsto x_ke(\beta,\gamma) & (1\le k\le m), \\
&1\tens x_k\mapsto x_{m+k}e(\beta,\gamma) & (1\le k\le n), \\
&\tau_k\tens 1\mapsto \tau_ke(\beta,\gamma) & (1\le k<m), \\
&1\tens \tau_k\mapsto \tau_{m+k}e(\beta,\gamma) & (1\le k<n).
\eneqn
Here $\mu*\nu$ is the concatenation of $\mu$ and $\nu$,
i.e., 
$$\mu*\nu=(\mu_1,\ldots,\mu_m,\nu_1,\ldots,\nu_n).$$

\medskip
For an $R(\beta)$-module $M$ and an $R(\gamma)$-module $N$,
we define their \emph{convolution product}
$M\conv N$ by
\eq&&
M\conv N=R(\beta + \gamma) e(\beta,\gamma)
\tens_{R(\beta )\otimes R( \gamma)}(M\otimes N). 
\eneq
Set $m=\height{\beta}$ and $n=\height{\gamma}$.
Set $$\sym_{m,n}\seteq\set{w\in\sym_{m+n}}{\text{\;$w\vert_{[1,m]}$ and
$w\vert_{[m+1,m+n]}$ are increasing}}.$$
Here $[a,b]\seteq\set{k\in\Z}{a\le k\le b}$.
Then we have
\eq
&&M\conv N=\soplus_{w\in\sym_{m,n}}\tau_w(M\tens N).
\label{eq:conv}
\eneq
We also have (see \cite[Theorem 2.2 (2)]{LV11})
\eq
&&\bl M\conv N\br^*\simeq N^*\conv M^*.
\label{eq: dual conv}
\eneq

\subsection{$R$-matrices for \KLRs}
\subsubsection{Intertwiners}\label{subsec:int}
 For $\height{\beta}=n$ and $1\le a<n$,  we define $\vphi_a\in R( \beta)$
by
\eq&&\ba{l}
  \vphi_a e(\nu)=
\begin{cases}
  \bl\tau_ax_a-x_{a}\tau_a\br e(\nu)\\
\hs{4ex}= \bl x_{a+1}\tau_a-\tau_ax_{a+1}\br e(\nu)\\
\hs{8ex}  =\bl\tau_a(x_a-x_{a+1})+1\br e(\nu)\\
\hs{12ex}=\bl(x_{a+1}-x_{a})\tau_a-1\br e(\nu)
 & \text{if $\nu_a=\nu_{a+1}$,} \\[2ex]
\tau_ae(\nu)& \text{otherwise.}
\end{cases}
\label{def:int} \ea \eneq
They are called the {\em intertwiners}.

\Lemma
\label{lem:ga} \hfill
\bnum
\item $\vphi_a^2e(\nu)=\bl
Q_{\nu_a,\nu_{a+1}}(x_a,x_{a+1})+\delta_{\nu_a,\nu_{a+1}}\br e(\nu)$.
\item
$\{\vphi_k\}_{1\le k<n}$ satisfies the braid relation.
\item
 For $w\in \sym_n$, let $w=s_{a_1}\cdots s_{a_\ell}$ be a reduced expression of $w$ and
set $\vphi_w=\vphi_{a_1}\cdots\vphi_{a_\ell}$.
Then $\vphi_w$ does not depend on the choice of reduced expressions of $w$.

\item For $w\in \sym_n$ and $1\le k\le n$, we have $\vphi_w x_k=x_{w(k)}\vphi_w$.
\item
For $w\in \sym_n$ and $1\le k<n$,
if $w(k+1)=w(k)+1$, then $\vphi_w\tau_k=\tau_{w(k)}\vphi_w$.\label{ga5}
\ee
\enlemma

For $m,n\in\Z_{\ge0}$,
let us denote by $w[{m,n}]$ the element of $\sym_{m+n}$  defined by
\eq
&&w[{m,n}](k)=\begin{cases}k+n&\text{if $1\le k\le m$,}\\
k-m&\text{if $m<k\le m+n$.}\end{cases}
\eneq

Let $\beta,\gamma\in \rtl^+$ with $\height{\beta}=m$, $\height{\gamma}=n$,
 and let
$M$ be an $R(\beta)$-module and $N$ an $R(\gamma)$-module. Then the
map $M\tens N\to N\conv M$ given by
$u\tens v\longmapsto \vphi_{w[n,m]}(v\tens u)$
is  $R( \beta)\tens R(\gamma )$-linear by the above lemma, and it extends
to an $R( \beta +\gamma)$-module  homomorphism 
\eq &&R_{M,N}\col
M\conv N\To N\conv M. \eneq
 Then we obtain the following commutative diagrams:
\eq&&\xymatrix@C=9ex@R=5ex{L\conv M\conv N\ar[r]^{R_{L,M}}\ar[dr]_{R_{L,M\circ N}}
&M\conv L\conv N\ar[d]^{R_{L,N}}\\
&M\conv N\conv L}\qtext{and}\quad
\xymatrix@C=9ex{L\conv M\conv N\ar[r]^{R_{M,N}}\ar[dr]_{R_{L\circ M,N}}
&L\conv N\conv M\ar[d]^{R_{L,N}}\\
&N\conv L\conv M\,.}
\label{eq:Rcom}
\eneq

\subsubsection{Spectral parameters}\label{subsec:spec}

\Def For $\beta\in\rtl^+$, the \KLR\ $R( \beta )$ is called {\em symmetric} if
$Q_{i,j}(u,v)$ is a polynomial in $u-v$ for all $i,j\in\supp(\beta)$.
Here, we set $\supp(\beta)=\set{i_k}{1\le k\le n}$
for $\beta=\sum_{k=1}^n\al_{i_k}$.
\edf

\medskip
Assume that the \KLR\ $R(\beta)$ is symmetric.
Let $z$ be an indeterminate, and
let $\psi_z$ be the algebra homomorphism
\eqn
&&\psi_z\col R( \beta )\to \cor[z]\tens R( \beta )
\eneqn
given by
$$\psi_z(x_k)=x_k+z,\quad\psi_z(\tau_k)=\tau_k, \quad\psi_z(e(\nu))=e(\nu).$$

For an $R( \beta )$-module $M$, we denote by $M_z$
the $\bl\cor[z]\tens R( \beta )\br$-module
$\cor[z]\tens M$ with the action of $R( \beta )$ twisted by $\psi_z$.
Namely,
\eq&&
\ba{l}
e(\nu)(a\tens u)=a\tens e(\nu)u,\\[1ex]
x_k(a\tens u)=(za)\tens u+a\tens (x_ku),\\[1ex]
\tau_k(a\tens u)=a\tens(\tau_k u)
\ea
\label{eq:spt1}
\eneq
for  $\nu\in I^\beta$,  $a\in \cor[z]$ and $u\in M$.
 For $u\in M$, we sometimes denote by $u_z$ the corresponding element $1\tens u$ of
the $R( \beta )$-module  $M_z$.

For a non-zero $M\in R(\beta)\smod$
and a non-zero $N\in R(\gamma)\smod$,
\eq&&\parbox{70ex}{%
let $s$ be the order of zeroes of $R_{M_z,N}\col  M_z\conv
N\To N\conv M_z$;
i.e., the largest non-negative integer such that the image of
$R_{M_z,N}$ is contained in $z^s(N\conv M_z)$.}
\label{def:s} \eneq
Note that such an $s$ exists because
$R_{M_z,N}$ does not vanish (\cite[Proposition 1.4.4 (iii)]{K^3}).

\Def Assume that $R(\beta)$ is symmetric.
For a non-zero $M\in R(\beta)\smod$ 
and a non-zero $N\in R(\gamma)\smod$,  
let $s$ be an integer as in \eqref{def:s}.
We define $$\rmat{M,N}\col M\conv N\to N\conv M$$
by  $$\rmat{M,N} = \bl z^{-s}R_{M_z,N}\br\vert_{z=0},$$
and call it the {\em renormalized $R$-matrix}. 
\edf
By the definition, the renormalized $R$-matrix
$\rmat{M,N}$ never vanishes.

We define also
$$\rmat{N,M}\col N\conv M\to M\conv N$$
by  $$\rmat{N,M} = \bl (-z)^{-t}R_{N,M_z}\br\vert_{z=0},$$
where $t$ is the multiplicity of zero of
$R_{N,M_z}$.

Note that if $R(\beta)$ and $R(\gamma)$ are symmetric, then
$s$ coincides with the multiplicity of zero of
$R_{M,N_z}$, and $\bl z^{-s}R_{M_z,N}\br\vert_{z=0}=\bl (-z)^{-s}R_{M,N_z}\br\vert_{z=0}$.

Indeed, we have
\eq
&&\ba{l}
R_{M_{z_1},N_{z_2}}\bl(u)_{z_1}\tens (v)_{z_2})\br
=\vphi_{w[n,m]}\bl (v)_{z_2}\tens(u)_{z_1})\br\\[1ex]
\hs{30ex}\in
\sum_{w,u',v'}\cor[z_1-z_2]\tau_w\bl (v')_{z_2}\tens (u')_{z_1}\br
\ea
\label{eq:-}\eneq
for $u\in M$ and $v\in N$.
Here $w$ ranges over
$$\sym_{n,m}\seteq\set{w\in\sym_{m+n}}{\text{$w\vert_{[1,n]}$ and $w\vert_{[n+1,n+m]}$
are strictly increasing}}$$
and $v'\in N$ and $u'\in M$. 
Hence, $\rmat{M,N}$ is well defined whenever
at least one of $R(\beta)$ and $R(\gamma)$ is symmetric.

The proof of \eqref{eq:-} will be given later in Section \ref{sec:Proof}.

\section{Quantum affine algebras}\label{sec:quantum}

In this section, we briefly review the representation theory of quantum affine algebras following \cite{AK, Kas02}.
When concerned with quantum affine algebras, 
we take the algebraic closure of $\C(q)$
in $\cup_{m >0}\C((q^{1/m}))$ as a base field $\cor$.
\subsection{Integrable modules}
Let $ I$ be an index set and $A=(a_{ij})_{i,j \in I}$ be a generalized
Cartan matrix of affine type.

 We choose $0\in I$ as the leftmost vertices in the tables
in \cite[{pages 54, 55}]{Kac} except $A^{(2)}_{2n}$-case in which
case we take the longest simple root as $\al_0$. Set $I_0
=I\setminus\{0\}$.

The weight lattice $P$ is given by
\begin{align*}
  P = \Bigl(\soplus_{i\in I}\Z \La_i\Bigr) \oplus \Z \delta,
\end{align*}
and the simple roots are given by
$$\al_i=\sum_{j\in I}a_{ji}\La_j+\delta(i=0)\delta.$$
The weight $\delta$ is called the imaginary root. 
There exist $d_i\in\Z_{>0}$ such that
$$\delta=\sum_{i\in I}d_i\al_i.$$
Note that $d_i=1$ for $i=0$.
The simple coroots $h_i\in P^\vee\seteq\Hom[\Z] (P,\Z)$ are given by
\eqn
&& \lan h_i,\La_j\ran=\delta_{ij},\quad \lan h_i,\delta\ran=0.
\eneqn
Hence we have $\lan h_i,\al_j\ran=a_{ij}$.

Let $c=\sum_{i\in I}c_i h_i$ be a unique element such that
$c_i\in\Z_{>0}$ and
$$\Z\, c=\set{h\in\soplus\nolimits_{i\in I}\Z h_i}{\text{$\lan h,\al_i\ran=0$ for any $i\in I$}}.$$
Let us take a $\Q$-valued symmetric bilinear form  $(\scbul,\scbul)$ 
on $P$
such that
$$\text{
$\lan h_i,\la\ran=\dfrac{2(\al_i,\la)}{(\al_i,\al_i)}$
and $(\delta, \la)=\lan c,\la\ran $ for any $\la\in P$.}
$$
Let $q$ be an indeterminate. For each $i \in I$, set $q_i = q^{(\al_i,\al_i)/2}$.

\begin{definition} \label{def:qgroup}
The {\em quantum group} $U_q(\g)$ associated with $(A,P)$ is the  
$\cor$-algebra  generated by $e_i,f_i$ $(i \in I)$ and
$q^{\la}$ $(\la \in P)$ satisfying following relations\,{\rm:}
\begin{equation*}
\begin{aligned}
& q^0=1,\ q^{\la} q^{\la'}=q^{\la+\la'} \ \ \text{for} \ \la,\la' \in P,\\
& q^{\la}e_i q^{-\la}= q^{(\la, \alpha_i)} e_i, \ \
          \ q^{\la}f_i q^{-\la} = q^{-(\la, \alpha_i)} f_i \ \ \text{for} \ \la \in P, 
i \in
          I, \\
& e_if_j - f_je_i = \delta_{ij} \dfrac{K_i -K^{-1}_i}{q_i- q^{-1}_i
}, \ \ \mbox{ where } K_i=q^{\al_i}, \\
& \sum^{1-a_{ij}}_{r=0} (-1)^r \left[\begin{matrix}1-a_{ij}
\\ r\\ \end{matrix} \right]_i e^{1-a_{ij}-r}_i
         e_j e^{r}_i =0 \quad \text{ if } i \ne j, \\
& \sum^{1-a_{ij}}_{r=0} (-1)^r \left[\begin{matrix}1-a_{ij}
\\ r\\ \end{matrix} \right]_i f^{1-a_{ij}-r}_if_j
        f^{r}_i=0 \quad \text{ if } i \ne j.
\end{aligned}
\end{equation*}
\end{definition}

Here, we set $[n]_i =\dfrac{ q^n_{i} - q^{-n}_{i} }{ q_{i} - q^{-1}_{i} },\quad
  [n]_i! = \prod^{n}_{k=1} [k]_i$ and
  $\left[\begin{matrix}m \\ n\\ \end{matrix} \right]_i= \dfrac{ [m]_i! }{[m-n]_i! [n]_i! }\;$
  for each $n \in \Z_{\ge 0}$, $i \in I$ and $m\ge n$.

 We denote by $\uqpg$ the subalgebra of $U_q(\mathfrak{g})$ generated by
$e_i$, $f_i$, $K_i^{\pm1}(i\in I)$, and call it {\em quantum affine algebra}.
The algebra $\uq$ has a Hopf algebra structure with the coproduct:
\begin{align}
\Delta(K_i)&=K_i\tens K_i, \nn\\
\Delta(e_i)&=e_i\tens K_i^{-1}+1\tens e_i,\\
\Delta(f_i)&=f_i\tens 1+K_i\tens f_i.\nn
\end{align}

Set
$$P_\cl=P/\Z\delta$$
and call it the {\em classical weight lattice}. Let $\cl\col P\to
P_\cl$ be the projection. Then $P_\cl=\soplus\nolimits_{i\in I}\Z\,\cl(\La_i)$.
Set $P^0_\cl = \set{\la\in P_\cl}{\lan c,\la\ran=0}\subset P_\cl$.

A $\uqpg$-module $M$ is called an {\em integrable module} if

\hs{5ex}\parbox[t]{70ex}{
 \bna
\item $M$ has a weight space decomposition
$$M = \bigoplus_{\lambda \in P_\cl} M_\lambda,$$
where  $M_{\lambda}= \set{ u \in M }{%
\text{$K_i u =q_i^{\lan h_i , \lambda \ran} u$ for all $i\in I$}}$,
\item  the actions of
 $e_i$ and $f_i$ on $M$ are locally nilpotent for any $i\in I$.
\ee}

\medskip
Let us denote by $\uqm$
the abelian tensor category of finite-dimensional integrable $\uqpg$-modules.

If $M$ is a  simple module in $\uqm$,
then there exists a non-zero vector $u\in M$
of weight $\la\in P_\cl^0$
such that $\la $ is dominant (i.e.,
$\lan h_i,\la\ran\ge0$ for any $i\in I_0$)
and all the weights of $M$ lie in
$\la-\sum_{i\in I_0}\Z_{\ge0}\al_i$.
We say that $\la$ is the {\em dominant extremal weight} of $M$ and $u$
is a {\em dominant extremal vector} of $M$. Note that
a dominant extremal vector of $M$ is unique up to a constant multiple.

Let $z$ be an indeterminate.
For a $\uqpg$-module $M$, let us denote by $M_z$
the module $\cor[z,z^{-1}]\tens M$ with the action of $\uqpg$ given by
$$e_i(u_z)=z^{\delta_{i,0}}(e_iu)_z, \quad
f_i(u_z)=z^{-\delta_{i,0}}(f_iu)_z, \quad
K_i(u_z)=(K_iu)_z.
$$
Here, for $u\in M$,  we denote by $u_z$ the element
$1\tens u\in \cor[z,z^{-1}]\tens M$.

\subsection{$R$-matrices}
 We recall the notion of $R$-matrices \cite[\S\;8]{Kas02}.
Let us choose the following {\em universal $R$-matrix}. 
Let us take a basis $\{P_\nu\}_\nu$ of $U_q^+(\g)$
and a basis $\{Q_\nu\}_\nu$ of $U_q^-(\g)$
dual to each other with respect to a suitable coupling between
$U_q^+(\g)$ and $U_q^-(\g)$.
Then for $\uq$-modules $M$ and $N$
define
\eq\label{def:univ}
R_{MN}^\univ(u\otimes v)=q^{(\wt(u),\wt(v))}
\sum_\nu P_\nu v\otimes Q_\nu u\,,
\eneq
so that
$R_{MN}^\univ$ gives a $\uq$-linear homomorphism
from $M\otimes N$ to $N\otimes M$
provided that the infinite sum has a meaning.

Let $M$ and $N$ be $\uq$-modules in $\uqm$, 
and let $z_1$ and $z_2$ be indeterminates. 
Then $R_{M_{z_1},N_{z_2}}^\univ$
converges in the $(z_2/z_1)$-adic topology.
Hence we obtain a morphism of
$\cor[[z_2/z_1]]\tens_{\cor[z_2/z_1]}\cor[z_1^{\pm1},z_2^{\pm1}]\tens\uq$-modules
$$R_{M_{z_1},N_{z_2}}^\univ\col \cor[[z_2/z_1]]\tens_{\cor[z_2/z_1]}
(M_{z_1}\tens N_{z_2})
\to \cor[[z_2/z_1]]\tens_{\cor[z_2/z_1]}
(N_{z_2}\tens M_{z_1}).
$$
If there exist $a\in\cor((z_2/z_1))$ and
 a $\cor[z_1^{\pm1},z_2^{\pm1}]\tens\uq$-linear homomorphism
$$R\col M_{z_1}\tens N_{z_2}\to N_{z_2}\tens M_{z_1}$$
such that $R^\univ_{M_{z_1},N_{z_2}}=a R$,
then we say that $R^\univ_{M_{z_1},N_{z_2}}$ is {\em rationally renormalizable.}

Now assume further that $M$ and $N$ are non-zero.
Then, we can choose $R$ so that, for any $c_1,c_2\in \cor^\times$,
the specialization of $R$ at $z_1=c_1$, $z_2=c_2$
$$R\vert_{z_1=c_1,z_2=c_2}\col M_{c_1}\tens N_{c_2}\to N_{c_2}\tens M_{c_1}$$
does not vanish.
Such an $R$ is unique up to a multiple of
$\cor[(z_1/z_2)^{\pm1}]^\times=\sqcup_{n\in\Z}\,\cor^\times z_1^nz_2^{-n}$.
We write
$$\rmat{M,N}\seteq R\vert_{z_1=z_2=1}\col M\tens N\to N\tens M,$$
and call it the {\em renormalized $R$-matrix}. 
The renormalized $R$-matrix $\rmat{M,N}$ 
is well defined up to a constant multiple
when $R^\univ_{M_{z_1},N_{z_2}}$ is rationally renormalizable.
By the definition, $\rmat{M,N}$ never vanishes.

\medskip
Now assume that $M_1$ and $M_2$ are simple $\uqpg$-modules in $\uqm$.
Then, the universal $R$-matrix 
$R^\univ_{(M_1)_{z_1},(M_2)_{z_2}}$ is rationally renormalizable.
More precisely, we have the following.
Let $u_1$ and $u_2$
be dominant extremal weight vectors of $M_1$ and $M_2$,
respectively. 
Then there exists $a(z_2/z_1)\in\cor[[z_2/z_1]]^\times$
such that
$$R^\univ_{(M_1)_{z_1},(M_2)_{z_2}}\bl(u_1)_{z_1}\tens (u_2)_{z_2}\br=
a(z_2/z_1)\bl(u_2)_{z_2}\tens(u_1)_{z_1}\br.$$
Then $\Rnorm_{M_1, M_2}\seteq a(z_2/z_1)^{-1}R^\univ_{(M_1)_{z_1},(M_2)_{z_2}}$
is a unique $\cor(z_1,z_2)\tens\uqpg$-module homomorphism
\begin{equation}\ba{l}
\Rnorm_{M_1, M_2} \col
\cor(z_1,z_2)\otimes_{\cor[z_1^{\pm1},z_2^{\pm1}]} \bl(M_1)_{z_1} \otimes (M_2)_{z_2}\br \\[2ex]
\hs{30ex}\To\cor(z_1,z_2)\otimes_{\cor[z_1^{\pm1},z_2^{\pm1}]}
 \bl(M_2)_{z_2} \otimes (M_1)_{z_1} \br
\ea\end{equation}
satisfying
\begin{equation}\Rnorm_{M_1, M_2}\bl(u_1)_{z_1} \otimes (u_2)_{z_2}\br = 
(u_2)_{z_2}\otimes (u_1)_{z_1}.
\end{equation}

Note that
$\cor(z_1,z_2)\otimes_{\cor[z_1^{\pm1},z_2^{\pm1}]}
 \big((M_1)_{z_1} \otimes (M_2)_{z_2} \big)$ is a simple
$\cor(z_1,z_2)\tens\uqpg$-module (\cite[Proposition 9.5]{Kas02}). 
We call $\Rnorm_{M_1, M_2}$ the
{\em normalized $R$-matrix}.

Let $d_{M_1,M_2}(u) \in \cor[u]$ be a monic polynomial of the
smallest degree such that the image of $d_{M_1,M_2}(z_2/z_1)
\Rnorm_{M_1, M_2}$ is contained in $(M_2)_{z_2} \otimes(M_1)_{z_1}$. We call $d_{M_1,M_2}(u)$ the {\em denominator of
$\Rnorm_{M_1, M_2}$}.
Then we have
\begin{equation}
d_{M_1,M_2}(z_2/z_1)\Rnorm_{M_1, M_2}
\col (M_1)_{z_1} \otimes (M_2)_{z_2} \To
(M_2)_{z_2} \otimes (M_1)_{z_1},
\end{equation}
and the renormalized $R$-matrix 
$$\rmat{M_1,M_2}\col M_1\tens M_2\To M_2\tens M_1$$
is equal to the specialization of
$d_{M_1,M_2}(z_2/z_1)\Rnorm_{M_1, M_2}$ at $z_1=z_2=1$
up to a constant multiple.

\smallskip

Note that $R^\univ$ satisfies the following properties:
the following diagrams commute
\eqn
\xymatrix@C=15ex{
M_1\tens M_2\tens N\ar[r]_{M_1\tens R^\univ_{M_2,N}}\ar@/^2pc/[rr]^
{R^\univ_{M_1\tens M_2,N}}
&M_1\tens N\tens M_2\ar[r]_{R^\univ_{M_1,N}\tens M_2}&
N\tens M_1\tens M_2,
}\\
\xymatrix@C=15ex{
M\tens N_1\tens N_2\ar[r]_{R^\univ_{M,N_1}\tens N_2}
\ar@/^2pc/[rr]^{R^\univ_{M, N_1\tens N_2}}
&N_1\tens M\tens N_2\ar[r]_{N_1\tens R^\univ_{M,N_2}}&
N_1\tens N_2\tens M
}
\eneqn
for $M$, $M_1$, $M_2$, $N$, $N_1$, $N_2$ in $\uqm$.
Hence, if $R^\univ_{(M_1)_{z_1},N_{z_2}}$ and $R^\univ_{(M_2)_{z_1},N_{z_2}}$
are rationally renormalizable, then
$R^\univ_{(M_1\tens M_2)_{z_1},N_{z_2}}$ is also rationally renormalizable.
Moreover, we have
\eq\text{
$\bl\rmat{M_1,N}\tens M_2\br\circ\bl M_1\tens \rmat{M_2,N}\br
=c\,\rmat{M_1\tens M_2,N}$ for some $c\in \cor$.}
\label{eq:rtens}
\eneq
Note that $c$ may vanish.

In particular, if $M_1$, $M_2$ and $N$ are simple modules in $\uqm$, then
$R^\univ_{(M_1\tens M_2)_{z_1},N_{z_2}}$ is  rationally renormalizable.

\section{Simple  heads and socles  of tensor products}
In this section we give a proof of Conjecture in Introduction
for the \KLR\ case and the quantum affine algebra case.

\subsection{Quiver Hecke algebra case}\label{sec:qHal}
We shall first discuss the \KLR\ case.

\Lemma\label{lem:sub}
Let $\beta_k\in\rtl^+$ and $M_k\in R(\beta_k)\gmod$
$(k=1,2,3)$.
Let $X$ be an $R(\beta_1+\beta_2)$-submodule of $M_1\conv M_2$ and
$Y$ an $R(\beta_2+\beta_3)$-submodule of $M_2\conv M_3$ 
such that $X\conv M_3\subset M_1\conv Y$ as submodules of
$M_1\conv M_2\conv M_3$.
Then there exists an $R(\beta_2)$-submodule $N$ of $M_2$
such that
$X\subset M_1\conv N$ and $N\conv M_3\subset Y$.
\enlemma

\Proof
Set $n_k=\height{\beta_k}$.
Set $N=\set{u\in M_2}{u\tens M_3\subset Y}$.
Then $N$ is the largest $R(\beta_2)$-submodule of $M_2$ such that
$N\conv M_3\subset Y$.
Let us show $X\subset M_1\conv N$.
Let us take a basis  $\{v_a\}_{a\in A}$ of $M_1$.

By \eqref{eq:conv}, we have
$$M_1\conv M_2=\soplus_{w\in\sym_{n_1,n_2}}\tau_w(M_1\tens M_2).$$
Hence, any $u\in X$ can be uniquely written as
$$u=\sum_{w\in\sym_{n_1,n_2},\,a\in A}\tau_w(v_a\tens u_{a,w})$$
with $u_{a,w}\in M_2$. Then, for any $s\in M_3$, we have 
$$u\tens s=\sum_{w\in\sym_{n_1,n_2},\,a\in A}\tau_w(v_a\tens u_{a,w}\tens s)
\in X\conv M_3\subset  M_1 \conv Y.$$ Since
$$M_1 \conv Y =\soplus_{w\in\sym_{n_1,n_2+n_3}}\tau_w(M_1 \tens Y)$$
and $\sym_{n_1,n_2}\subset\sym_{n_1,n_2+n_3}$,
we have
$$u_{a,w}\tens s\in Y \quad
\text{for any $a\in A$ and $w\in\sym_{n_1,n_2}$.}$$
Therefore we have
$u_{a,w}\in N$.
\QED

\Th\label{th:main}
Let $\beta,\gamma\in \rtl^+$
and $M\in R(\beta)\smod$ and $N\in R(\gamma)\smod$. 
 We assume further the following condition:
\eq&&\parbox{70ex}{\bna
\item $R(\beta)$ is symmetric and
$\rmat{M,M}\in\cor \id_{M\circ M}$,
\item $M$ is non-zero,
\item $N$ is a simple $R(\gamma)$-module.
\ee}
\label{cond:real}
\eneq
Then we have
\bnum
\item
$M\conv N$ has a simple socle and a simple head.
Similarly, $N\conv M$ has a simple socle and a simple head.
\item Moreover,
$\Im (\rmat{N, M})$ is equal to the socle of $M\conv N$ and
also equal to the head of $ N\conv Mx$.
Similarly, $\Im (\rmat{M, N})$ is equal to the socle of $N \conv M$ and
to the head of $M\conv N $.
\ee
In particular, $M$ is a simple module.
\enth
 \Proof
Let us show that 
$\Im (\rmat{N , M})$ is a unique simple submodule of
$M\conv N $. Let $S\subset M\conv N $ be an arbitrary
non-zero $R(\beta+\gamma)$-submodule. Let $m$ and $m'$ 
be the multiplicity of zero of $R_{N ,(M)_z}\col  N \conv
(M)_z\to (M)_z\conv N $ and $R_{M,(M)_z}\col  M\conv
(M)_z\to (M)_z\conv M$ at $z=0$, respectively. Then by the definition,
$\rmat{N ,M}=(z^{-m}R_{N ,(M)_z})\vert_{z=0}\col N \conv M\to
M\conv N $ and
$\rmat{M,M}=(z^{-m'}R_{M,(M)_z})\vert_{z=0}\col M\conv
M\to M\conv M$. 
Now, we
have a commutative diagram

\eqn&&\hs{-4.5ex}\xymatrix@C=15ex@R=5ex{
\wb{S\conv (M)_z}\ar[rr]^{z^{-m-m'}R_{S, (M)_z}}\xmono[d]&&
\wb{(M)_z\conv S}\xmono[d]\\
M\conv N \conv (M)_z\ar[r]^{M\circ z^{-m}R_{N ,(M)_z}}
&M\conv (M)_z\conv N \ar[r]^{z^{-m'}R_{M,(M)_z}\circ N }
&(M)_z\conv M\conv N .
}
\eneqn
Therefore
$z^{-m-m'}R_{S, (M)_z}\col S\conv (M)_z\to (M)_z\conv S$ is well-defined,
and we obtain the following commutative diagram
by specializing the above diagram at $z=0$.
$$
\xymatrix@C=10ex@R=5ex{
\wb{S\conv M}\ar[rr]\xmono[d]&&\wb{M\conv S}\xmono[d]\\
M\conv N  \conv M\ar[r]^{M\circ  \rmats{N ,M}}
&M\conv M\conv N \ar[r]^{\id_{M\circ M\circ N }}
&M
\conv M\conv N \,.
}$$
Here, we have used the assumption
that $\rmat{M,M}$ is equal to $\id_{M\circ M}$ up to a constant multiple.

\noi
Hence we obtain
$\bl M\conv \rmat{N ,M}\br(S\conv M)\subset M\conv S$,
or equivalently
$$S\conv M\subset M\conv (\rmat{N ,M})^{-1}(S).$$
By the preceding lemma, there exists an
$R(\gamma)$-submodule $K$ of $N $ such that
$S\subset M\conv K$ and $K\conv M\subset (\rmat{N ,M})^{-1}(S)$.
By the first inclusion, we have $K\not=0$.
Since $N $ is simple, we have
$K=N $ and we obtain $N \conv M\subset (\rmat{N ,M})^{-1}(S)$,
or equivalently,
$\Im (\rmat{N ,M})\subset S$.
Noting that $S$ is an arbitrary non-zero submodule of $M\conv N $,
we conclude that $\Im (\rmat{N ,M})$ is a unique simple submodule of
$M\conv N $.

The proof of the other statements in (i) and (ii) is similar.

The simplicity of $M$ follows from (i) and (ii) by taking 
the one-dimensional $R(0)$-module $\cor$ as $N $. Note that
$\rmat{M,\cor}$ and $\rmat{\cor,M}$ coincide with
the identity morphism $\id_M$.
\QED

A simple $R(\beta)$-module $M$ is called {\em real}
if $M\conv M$ is simple

Then the following corollary is an immediate consequence of 
Theorem~\ref{th:main}.

\Cor Assume that $R(\beta)$ is symmetric
and $M$ is a non-zero $R(\beta)$-module in $R(\beta)\smod$. 
Then the following conditions are equivalent:
\bna
\item $M$ is a real simple $R(\beta)$-module,
\item$\rmat{M,M}\in\cor\id_{M\circ M}$,
\item $\End_{R(2\beta)}(M\conv M)\simeq\cor\id_{M\circ M}$.
\ee
\encor

We have also the following corollary.
\Cor
If $R(\beta)$ is symmetric and $M$ is a real simple $R(\beta)$-module, 
then $M^{\circ\, n}\seteq\overbrace{M\conv\cdots \conv M}^n$ is 
a simple $R(n\beta)$-module for any $n\ge 1$.
\encor 
\Proof The quiver Hecke algebra version of \eqref{eq:rtens}
implies that $\rmat{M^{\circ\, m},M^{\circ\, n}}$ is equal to $\id_{M^{\circ (m+n)}}$ 
up to a constant multiple.
\QED

Thus we have established the first statement of Conjecture in the introduction
in the quiver Hecke algebra case.

\Lemma\label{lem:XY}
Let $\beta,\gamma\in\rtl^+$, and
let $M\in R(\beta)\smod$ and $L\in R(\beta+\gamma)\smod$.
Then there exist $X$, $Y\in R(\gamma)\smod$
satisfying the following universal properties\,{\rm:}
\eq&&
\Hom[{R(\beta+\gamma)}](M\conv Z,L)\simeq\Hom[{R(\gamma)}](Z,X),
\label{eq:X}\\
&&\Hom[{R(\beta+\gamma)}](L,Z\conv M)\simeq\Hom[{R(\gamma)}](Y, Z)\quad
\label{eq:Y}
\eneq
functorially in $Z\in R(\gamma)\smod$.
\enlemma
\Proof
Set $X=\Hom[{R(\beta+\gamma)}](M\conv R(\gamma),L)$. Then we have
\begin{align*}
\Hom[{R(\beta+\gamma)}](M\conv Z,L)
&\simeq\Hom[{R(\beta)\tens R(\gamma)}](M\tens Z,L)\\
&\simeq\Hom[{R(\gamma)}]\bl Z,\Hom[{R(\beta)}](M,L)\br.
\end{align*}
Similarly set $Y=\Bigl(\Hom[{R(\beta+\gamma)}](M^*\conv R(\gamma),L^*)\Bigr)^*$.
Then we have by using \eqref{eq: dual conv}
\begin{align*}
\Hom[{R(\beta+\gamma)}](L, Z\conv M)
&\simeq\Hom[{R(\beta+\gamma)}](M^*\conv Z^*,L^*)\\
&\simeq\Hom[{R(\beta)\tens R(\gamma)}](M^*\tens Z^*,L^*)\\
&\simeq\Hom[{R(\gamma)}](Z^*,Y^*)
\simeq\Hom[{R(\gamma)}](Y,Z).
\end{align*}
\QED

\Prop Let $\beta,\gamma\in\rtl^+$. Assume that $R(\beta)$ is symmetric,
and let $M$ be a real simple module in $R(\beta)\smod$, 
and $L$ a simple module in $R(\beta+\gamma)\smod$. 
Then the $R(\gamma)$-module $X\seteq\Hom[{R(\beta+\gamma)}](M\conv R(\gamma),L)$
is either zero or
has a simple socle.
\enprop
\Proof
The $R(\gamma)$-module $X$ satisfies the functorial property
\eqref{eq:X}.
Assume that $X\not=0$.
Let $p\col M\conv X\to L$ be the canonical morphism.
Since $L$ is simple, it is an epimorphism.
Let $Y$ be as in Lemma~\ref{lem:XY},
and let $i\col L\to Y\circ M$ be the canonical morphism. 
For an arbitrary simple $R(\gamma)$-submodule $S$ of $X$,
since $\Hom[{R(\beta+\gamma)}](M\conv S,L)\simeq\Hom[{R(\gamma)}](S,X)$,
the composition $M\conv S\to M\conv X\To[p]L$ does not vanish.
 Hence, by Theorem~\ref{th:main},
  $L$ is  the simple head of $M\conv S$ and
is the simple socle of $S\conv M$. Moreover, 
$L \cong \Im(\rmat{M,S})$.
Since the monomorphism $L\to S\conv M$ factors through $i$
by \eqref{eq:Y},
the morphism $i\col L\to Y\circ M$ is a monomorphism.

As in the proof of Theorem~\ref{th:main},
we have a commutative diagram
$$\xymatrix@C=15ex@R=5ex{
\wb{M\conv L}\ar[r]\xmono[d]^{M\circ i}&\wb{L\conv M}\xmono[d]^{i\circ M}\\
M\conv Y\conv M\ar[r]^{\rmats{M,Y}\circ M}&Y\conv M\conv M\,.
}
$$
Then we obtain $M\conv i(L)\subset (\rmat{M,Y})^{-1}(i(L))\conv M$.
Hence, by Lemma~\ref{lem:sub},
there exists an $R(\gamma)$-submodule $Z$ of $Y$ such that
$\rmat{M,Y}(M\conv Z)\subset i(L)$
and $i(L) \subset Z\conv M$.
The last inclusion induces a morphism
$L\to Z\conv M$ and it induces a morphism $Y\to Z$ by \eqref{eq:Y}. 
Since the composition $Y\to Z\to Y$ is the identity again by \eqref{eq:Y}, 
we have $Z=Y$. Hence $\Im(\rmat{M,Y})\subset i(L)$, which gives
the commutative diagram 
$$\xymatrix@C=12ex{
M\conv Y\ar[r]\ar@/^1.8pc/[rr]^{\rmat{M,Y}}&L\;\xmono[r]_{i}&Y\conv M\,.
}
$$
By the argument dual to the above one
(also see the proof of Proposition~\ref{prop:uniqueend}), 
we have a commutative diagram
$$\xymatrix@C=12ex{
M\conv X\xepi[r]_p\ar[r]\ar@/^1.8pc/[rr]^{\rmat{M,X}}&L\;\ar[r]_\xi&X\conv M\,.
}
$$
Hence $\xi\col L\to X\conv M$ is a monomorphism,
and $\Im(\rmat{M,X})$ is isomorphic to $L$. 
By \eqref{eq:Y}, there exists a unique morphism $\vphi\col Y\to X$ such that
$\xi$ factors as 
$$\xymatrix@C=12ex
{\mbox{$L\;\;$}\xmono[r]_-i\ar@/^1.8pc/[rr]^-\xi
&Y\conv M\;\ar[r]_-{\vphi\circ M}&X\conv M.
}$$ 
Let us show that $\Im(\vphi)$ is a unique simple submodule of $X$.
In order to see this, 
let $S$ be an arbitrary simple $R(\beta)$-submodule of $X$.
We have seen that $L$ is isomorphic to the head of $M\conv S$ and 
isomorphic to $\Im(\rmat{M,S})$. 
Since the composition $M\conv S\to M\conv X\To[{\rmat{M,X}}]X\conv M$
does not vanish, 
we have a commutative diagram by \cite[Lemma~1.4.8]{K^3}
$$\xymatrix@C=12ex@R=4ex
{\wb{M\conv S}\ar[r]^{\rmat{M,S}}\xmono[d]&\wb{S\conv M}\xmono[d]\\
M\conv X\ar[r]^{\rmat{M,X}}&X\conv M.}
$$
Since $\Im(\rmat{M,S})\simeq\Im(\rmat{M,X})\simeq L$,
the morphism $\xi\col L\to X\circ M$ factors as 
$L\to S\circ M\to X\circ M$.
Hence \eqref{eq:Y} implies that
$\vphi\col Y\to X$ factors through $Y\to S\to X$.
Thus we obtain $\Im(\vphi)\subset S$.
Since $S$ is an arbitrary simple submodule of $X$,
we conclude that
$\Im(\vphi)$ is a unique simple submodule of $X$.
\QED

Let $\beta,\gamma\in\rtl^+$. For
a simple $R(\beta)$-module $M$ and a simple $R(\gamma)$-module $N$, 
let us denote by $M\hconv N$
the head of $M\conv N$. 
\Cor\label{cor:inj}  Let $\beta,\gamma\in\rtl^+$. 
Assume that $R(\beta)$ is symmetric, 
and let $M$ be a real simple module in $R(\beta)\smod$.
Then, the map $N\mapsto M\hconv N$ is injective
from the set of the isomorphism classes of simple objects of
$R(\gamma)\smod$ 
to the set of the isomorphism classes of  simple objects of 
$R(\beta+\gamma)\smod$.
\encor
\Proof
Indeed, for a simple $R(\gamma)$-module $N$, 
$M\hconv N$ is a simple $R(\beta+\gamma)$-module
by Theorem~\ref{th:main}, and 
$N\subset X\seteq\Hom[{R(\beta+\gamma)}](M\conv R(\gamma),M\hconv N)$
is the socle of $X$ by the preceding proposition.
\QED
If $L(i)$ is the one-dimensional simple $R(\al_i)$-module, then
$L(i)$ is real and
$M\hconv L(i)$ corresponds to the crystal operator $\tilde{f}_iM$
and $L(i)\hconv M$ to the dual crystal operator $\tilde{f_i}^{\vee} M$
in \cite{LV11}. Hence, $\hconv$ is a generalization of the crystal operator as suggested in \cite{L03}.

\Prop\label{prop:uniqueend}
Let $\beta,\gamma\in\rtl^+$.
Assume that $R(\beta)$ is symmetric,
and let $M$ be a real simple module in $R(\beta)\smod$,
and $N$ a simple module in $R(\gamma)\smod$.
Then we have $\End_{R(\beta+\gamma)}(M\conv N)\simeq\cor\id_{M\circ N}$.
\enprop
\Proof 
Set $L=M\conv N$. 
Let $X,Y\in R(\gamma)\smod$ be as in Lemma~\ref{lem:XY}.
Let $p\col M\conv X\to L$ and $i\col L\to Y\conv M$ be the canonical morphisms.
Then the isomorphism
$M\conv N\to L$ induces a morphism $j\col N\to X$
such that the composition $M\conv N\To[M\circ j]M\conv X\To[p] L$
is that isomorphism.
Hence $p\col M\conv X\to L$ is an epimorphism.
Since $N$ is simple and $j$ does not vanish, the morphism 
$j\col N\to X$ is a monomorphism.

We have a commutative diagram
$$
\xymatrix@C=13ex@R=5ex{
\wb{M\conv M\conv X}\ar[r]^{\rmats{M,M}\circ X}\xepi[d]^{M\circ p}
&M\conv M\conv X\ar[r]^{M\circ\,\rmats{M,X}}&\wb{M\conv X\conv M}\xepi[d]^{p\circ M}\\
M\conv L\ar[rr]&&L\conv M.
}$$
Since $\rmat{M,M}$ is $\id_{M\circ M}$ up to a constant multiple,
we obtain a commutative diagram: 
$$
\xymatrix@C=13ex@R=5ex{
\wb{M\conv (M\conv X)}\ar[r]^{M\circ\,\rmats{M,X}}\xepi[d]^{M\circ p}
&\wb{M\conv  X\conv M}\xepi[d]^{p\circ M}\\
M\conv L\ar[r]
&L\conv M.
}$$
Therefore we have
$$M\circ\bl\rmat{M,X}(\Ker p)\br\subset (\Ker p)\conv M.$$
Hence Lemma~\ref{lem:sub} implies that there exists $Z\subset X$ such that
$\rmat{M,X}(\Ker p)\subset Z\conv M$ and
$M\conv Z\subset \Ker p$.
The last inclusion shows that
$M\conv Z\to M\conv X\to L$ vanishes.
Hence by \eqref{eq:X}, the morphism $Z\to X$ vanishes, or equivalently,
$Z=0$.
Hence we have $\rmat{M,X}(\Ker p)=0$.
Therefore $\rmat{M,X}$ factors through $p$:
$$\xymatrix@C=12ex{
M\conv X\xepi[r]_-p\ar[r]\ar@/^1.8pc/[rr]^{\rmat{M,X}}&L\;\ar[r]_-{\xi}&X\conv M\,.
}
$$
Since $\rmat{M,X} \neq 0$, the morphism $\xi$ does not vanish. By \eqref{eq:Y}, there exists
$\vphi\col Y\to X$ such that
$\xi\col L\to X\conv M$ coincides with the composition
$L\To[\;i\;]Y\conv M\To[\vphi\,\circ\, M]X\circ M$.
Then we have a commutative diagram with the solid arrows:
$$\xymatrix@C=10ex@R=1.5ex{
\wb{M\conv N}\ar[rr]^{\rmat{M,N}}\xmono[dd]^{M\circ j}
\ar@{..>}[dr]^(.62){\text{$\sim$}}
&&\wb{N\conv M}\xmono[dd]^{j\circ M}\\
&L\ar@{..>}[dr]^-\xi\\
\wb{M\conv X}\ar[rr]_{\rmat{M,X}}\ar@{..>>}[ur]^-{p}&&X\conv M.
}$$
Indeed, the commutativity follows from \cite[Lemma~1.4.8]{K^3}
and the fact that
the composition $M\conv N\To[M\circ j]M\conv X\To[{\rmat{M,X}}
]X\conv M$ does not vanish because 
it coincides with
$M\conv N\isoto L\To[\xi]X\conv M$.

Thus $\xi\col L\to X\conv M$
coincides with
the composition
\eqn&&
L\simeq M\conv N\To[\;\rmat{M,N}\;]N\conv M\To[\ j\circ M\ ] X\conv M.
\label{eq:YNX}
\eneqn
Hence \eqref{eq:Y} implies that $\vphi\col Y\to X$ decomposes as 
$$\xymatrix{Y\ar[r]^-{\psi}& N\;\xmono[r]^-{\;j\;} & X.}$$
Since $N$ is simple, $\psi$ is an epimorphism,
and we conclude that $N$ is the image of $\vphi\col Y\to X$.

\medskip
Now let us prove that any $f\in\End_{R(\beta+\gamma)}(L)$
satisfies $f\in\cor\id_{L}$. 
By the universal properties \eqref{eq:X} and \eqref{eq:Y},
the endomorphism $f$ induces endomorphisms $f_X\in\End_{R(\gamma)}(X)$ 
and $f_Y\in\End_{R(\gamma)}(Y)$ such that the diagrams 
with the solid arrows 
\eq&&
\ba{l}\xymatrix@C=10ex@R=5ex{\wb{{M\conv X}}\xepi[r]_-p\ar[d]^{M\circ f_X}
\ar@{..>}@/^1.5pc/[rr]|{\rmat{M,X}}
& L\ar[d]^-{f}\ar@{..>}[r]_-\xi&X\conv M\ar[d]^{f_X\circ M}\\
{M\conv X}\xepi[r]^-p\ar@{..>}@/_1.5pc/[rr]|{\rmat{M,X}}
&L\ar@{..>}[r]^-\xi&X\conv M}\ea
\quad\text{and}\quad
\ba{l}\xymatrix@C=10ex@R=5ex
{L\ar[r]^-{i} \ar[d]^{f}&Y\conv M\ar[d]^-{f_Y\circ M}\\
L\ar[r]^-{i}&Y\conv M}\ea
\label{dia:fXY}
\eneq
commute.
Since $\rmat{M,X} $ commutes with $f$,
 the left diagram with dotted arrows commutes.
Hence, the following diagram with the solid arrows
\eq&&\ba{l}\xymatrix@C=10ex@R=5ex
{Y\xepi[r]_-{\psi}\ar[d]^{f_Y}\ar@/^1.3pc/[rr]^-{\vphi}
&N\;\ar@{..>}[d]^-{f_N}\xmono[r]_-{j}&X\ar[d]^{f_X}\\
Y\xepi[r]^-{\psi}\ar@/_1.3pc/[rr]_-{\vphi}&N\;\xmono[r]^-{j}&X
}\ea\label{dia:fphi}
\eneq 
commutes.
Then we can add the dotted arrow $f_N$ 
so that the whole diagram \eqref{dia:fphi} commutes.
Since $N$ is simple, we have
$f_N=c\id_N$ for some $c\in\cor$.
By replacing $f$ with $f-c\,\id_L$,
we may assume from the beginning that $f_N=0$.
Then $f_X\circ j=0$. Now, $f=0$ follows from the commutativity of the diagram
$$
\xymatrix@C=10ex@R=5ex{
M\conv N\ar@/^2pc/[rr]^-{\mbox{$\sim$}}
\ar[r]^-{M\circ j}\ar@{..>}[dr]_-{0}&M\conv X\ar[r]^-p\ar[d]^{M\circ f_X}
& L\ar[d]^-{f}\\
&M\conv X\ar[r]^-p&L.}
$$
\QED

\Cor\label{cor:main} 
Let $\beta,\gamma\in\rtl^+$, and assume that $R(\beta)$ is symmetric.
Let $M$ be a real simple module in $R(\beta)\smod$, 
and $N$ a simple module in $R(\gamma)\smod$. 
\bnum
\item
If the head of $M\conv N$ and the socle of $M\conv N$ are isomorphic,
then
$M\conv N$ is simple and $M\conv N\simeq N\conv M$. 
\item If $M\conv N\simeq N\conv M$, then
$M\conv N$ is simple. Conversely, if $M\conv N$ is simple,
then $M\conv N\simeq N\conv M$.
\ee
\encor
\Proof
(i) Let $S$ be the head of $M\conv N$ and the socle of $M\conv N$.
Then $S$ is simple.
Now we have the morphisms
$$M\conv N\epito S\monoto M\conv N.$$
By the preceding proposition, the composition is equal to $\id_{M\circ N}$
up to a constant multiple.
Hence $M\conv N$ and $N\conv M$ are isomorphic to $S$.  

\noi
(ii) Assume first $M\conv N\simeq N\conv M$.
Then the simplicity of $M\conv N$ immediately follows from (i) because 
the socle of $M\conv N$ is isomorphic to
the head of $N\conv M$ by Theorem~\ref{th:main}. 

If $M\conv N$ is simple, then $\rmat{M,N}$ is injective.
Since $\dim(M\conv N)=\dim(N\conv M)$,
$\rmat{M,N}\col M\conv N\to N\conv M$ is an isomorphism.
\QED

Note that, when $R(\beta)$ and $R(\gamma)$ are symmetric,
for a real simple $R(\beta)$-module $M$  and
a real simple $R(\gamma)$-module $N$, their convolution
$M\conv N$ is real simple if $M\conv N\simeq N\conv M$.

\subsection{Quantum affine algebra  case}

The similar results to Theorem~\ref{th:main},
Corollary~\ref{cor:inj} and Corollary~\ref{cor:main} hold also for
quantum affine algebras.
Let $\uqpg$ be the quantum affine algebra as in \S\;\ref{sec:quantum}.
Recall that $\uqm$ denotes the category of 
finite-dimensional integrable $\uq$-modules.

First note that
the following lemma, an analogue of Lemma~\ref{lem:sub}
in the quantum affine algebra case, is almost trivial.
Indeed, the similar result holds for any rigid monoidal category
which is abelian and the tensor functor is additive.

\Lemma Let $M_k$ be a module
in $\uqm$ $(k=1,2,3)$. 
Let $X$ be a $\uqpg$-submodule of $M_1\tens M_2$ and
$Y$ a $\uqpg$-submodule of $M_2\tens M_3$ 
such that $X\tens M_3\subset M_1\tens Y$ as submodules of
$M_1\tens M_2\tens M_3$.
Then there exists a $\uqpg$-submodule $N$ of $M_2$
such that
$X\subset M_1\tens N$ and $N\tens M_3\subset Y$.
\enlemma

\Cor\label{cor:rcomp}\hfill
\bnum
\item Let $M_k$ be a module in $\uqm$ $(k=1,2,3)$, and
let $\vphi_1\col L\to M_1\tens M_2$
and $\vphi_2\col M_2\tens M_3\to L'$
be {\em non-zero} morphisms.
Assume further that $M_2$ is a simple module.
Then the composition
\eq&&
L\tens M_3\To[\vphi_1\tens M_3] M_1\tens M_2\tens M_3\To[M_1\tens\vphi_2]
M_1\tens L'
\label{eq:comp}
\eneq
does not vanish.
\item
Let $M$, $N_1$ and $N_2$ be simple modules in $\uqm$. 
Then the following diagram commutes up to a constant multiple{\rm:} 
$$\xymatrix@C=12ex
{M\tens N_1\tens N_2\ar[r]_{\rmats{M,N_1}\tens N_2}\ar@/^2pc/[rr]|-{\rmats{M,\,N_1\tens N_2}}
&N_1\tens M\tens N_2\ar[r]_{N_1\tens\,\rmats{M,N_2}}&N_1\tens N_2\tens M.
}$$ 
\ee
\encor
\Proof
(i)\ Assume that the composition \eqref{eq:comp} vanishes.
Then we have $\Im\vphi_1\tens M_3\subset M_1\tens \Ker\vphi_2$.
Hence, by the preceding lemma, there exists
$N\subset M_2$ such that
$\Im\vphi_1\subset M_1\tens N$ and
$N\tens M_3\subset \Ker\vphi_2$.
The first inclusion implies $N\not=0$ and the last inclusion implies $N\not=M_2$.
It contradicts the simplicity of $M_2$.

\smallskip
\noi
(ii)\ By (i) $(N_1\tens\rmat{M,N_2})\circ(\rmat{M,N_1}\tens N_2)$
does not vanish.
Hence it is equal to $\rmat{M,\,N_1\tens N_2}$ up to a constant multiple
by \eqref{eq:rtens}.
\QED
Since the proof of the following theorem is similar
to the \KLR\ case, we just state the result
omitting its proof.

\Th\label{th:main A}
Let $M$ and $N$ be simple modules in $\uqm$. 
We assume further
\eq&&\text{$\rmat{M,M}\in\cor \id_{M\tens M}$.}
\label{cond:real1}
\eneq
Then we have
\bnum
\item
$M\tens N$ has a simple socle and a simple head.
\item Moreover,
$\Im (\rmat{M, N})$ is equal to the head of $M\tens N$ and
is also equal to the socle of $N\tens M$.
\ee
\enth

Recall that a simple $\uq$-module $M$ is called {\em real} \/if $M\tens M$
is simple. Hence $M$ in Theorem~\ref{th:main A} is real. 

For a module $M$ in $\uqm$, let us denote by $\rM$
 and $\lM$  the right dual and the left dual of $M$, respectively.
Hence we have isomorphisms
\eq
&&\ba{l}
\Hom[\uq](M\tens X,Y)\simeq \Hom[\uq](X, \rM\tens Y),\\[1ex]
\Hom[\uq](X\tens \rM,Y)\simeq \Hom[\uq](X, Y\tens M),\\[1ex]
\Hom[\uq](\lM\tens X,Y)\simeq \Hom[\uq](X, M\tens Y),\\[1ex]
\Hom[\uq](X\tens M,Y)\simeq \Hom[\uq](X, Y\tens \lM)
\ea\label{eq:dual}
\eneq
functorial in $X$, $Y\in\uqm$.

\Cor Under the assumption of the theorem above, the head of\/
$\Im\rmat{M,N}\tens \rM$ is isomorphic to $N$. \encor
\Proof Set $S=\Im \rmat{M,N}$. Since $\Hom[\uq](S, N\tens M)\simeq
\Hom[\uq](S\tens \rM,N)$, there exists a non-trivial morphism
$S\tens \rM\to N$. Since $N$ is simple, we have an
epimorphism
$$S\tens \rM\epito N.$$
Since $\rM\tens \rM\simeq \!{}^*(M\tens M)$ is a simple module,
the tensor product $S\tens \rM$ has a simple
head by the preceding theorem. Hence, we obtain the desired result.
\QED

For simple $\uq$-modules $M$ and $N$, let us denote by $M\hconv N$
the head of $M\tens N$. 
\Cor\label{cor:bj} Let $M$ be a real simple module in $\uqm$.
Then, the map $N\mapsto M\hconv N$ is bijective on the set of
the isomorphism classes of simple $\uqpg$-modules in $\uqm$, and its inverse is
given by $N\mapsto N\hconv \rM$. \encor

\Lemma Let $M$ be a real simple module in $\uqm$ and $N$
a simple module in $\uqm$. Then we have
$\End_{\,\uq}(M\tens N)\simeq\cor\id_{M\tens N}$.
\enlemma
\Proof
By Corollary~\ref{cor:rcomp},
we have a commutative diagram up to  a constant multiple
$$\xymatrix@C=12ex
{\lM\tens M\tens N\ar[r]_{\rmats{\lM,M}\tens N}\ar@/^2pc/[rr]|-{\rmat{\lM,M\tens N}}
&M\tens \lM\tens N\ar[r]_{M\tens \rmats{\lM,N}}&M\tens N\tens \lM.
}$$ 
By Theorem~\ref{th:main A},
$\Im(\rmat{\;\lM,M})$ is the simple socle of  $M \tens \lM$, and hence
$\rmat{\;\lM,M}$ is equal to the composition
$$\lM\tens M\To[\ \eps\ ]\one\To M\tens\lM$$
up to a constant multiple. 
Here $\one$ denotes the trivial representation of $\uq$.
Hence
we have a commutative diagram up to  a constant multiple
$$\xymatrix@C=12ex
{\lM\tens M\tens N\xepi[r]_{\eps\tens N}\ar@/^2pc/[rr]|-{\rmat{\;\lM,M\tens N}}
&\;N\;\xmono[r]&\;M\tens N\tens \lM.
}$$ 
Let $f\in \End_{\,\uq}(M\tens N)$.
Let us show that $f\in\cor\id_{M\tens N}$.
Since $\rmat{\;\lM,M\tens N}$ commutes with $f$, 
the following diagram with the solid arrows
\eq&&\ba{c}\xymatrix@C=12ex{
\lM\tens M\tens N\xepi[r]\ar[d]^-{\lM\tens f}&\;N\;\xmono[r]
\ar@{.>}[d]^-{f_N}
&\;M\tens N\tens \lM\ar[d]^-{f\tens \lM}\\
\lM\tens M\tens N\;\xepi[r]^-{\eps\tens N}&\;N\;\xmono[r]&\;M\tens N\tens \lM
}\ea\label{dia:MN}\eneq
is commutative.
Hence  we can add the dotted arrow $f_N$ 
so that the whole diagram \eqref{dia:MN} commutes. 
Since $N$ is simple, we have $f_N=c\id_N$ for some $c\in\cor$.
Then by replacing $f$ with $f-c\id_{M\tens N}$,
we may assume from the beginning that $f_N=0$.
Hence the composition
$$\lM\tens M\tens N\To[\lM\tens f]\lM\tens M\tens N
\To[\eps\tens N]N$$
vanishes.
Therefore \eqref{eq:dual} implies that
$M\tens N\To[\ f\ ]M\tens N$
vanishes.
\QED

\Cor\label{cor:mainQ} Let $M$ be a real simple module in $\uqm$, 
and $N$ a simple module in $\uqm$. 
\bnum
\item
If the head of $M\tens N$ and the socle of $M\tens N$ are isomorphic,
then
$M\tens N$ is simple and $M\tens N\simeq N\tens M$. 
\item If $M\tens N\simeq N\tens M$, then $M\tens N$ is simple.
\ee
\encor
This corollary follows from
the preceding lemma by an argument similar to the one in
the proof of Corollary~\ref{cor:main}.

\section{Proof of \eqref{eq:-}}\label{sec:Proof}
We shall show \eqref{eq:-}.
We keep the notations in Section\;\ref{sec:KLR}. 
We set
$$
\tx_{a,b}=\hs{-2ex}\sum_{\substack{\nu\in I^{\beta+\gamma},\\
\nu_a,\nu_b\in\supp(\beta)\cap\supp(\gamma)}}(x_a-x_b)e(\nu)\quad\text{and}
\quad \ttau_{c}=\hs{-2ex}\sum_{\substack{\nu\in I^{\beta+\gamma},\\
\nu_c\in\supp(\gamma),\,\nu_{c+1}\in\supp(\beta)}}\tau_ce(\nu)$$
for $1\le a,b\le m+n$ and $1\le c<m+n$.
They are elements of $R(\beta+\gamma)$.

We denote by $A$ the commutative subalgebra of $R(\beta+\gamma)$
generated by $\tx_{a,b}$ and $e(\nu)$ where
$1\le a<b\le m+n$ and $\nu\in I^{\beta+\gamma}$.
Let us denote by $\tR$ the subalgebra of $R(\beta+\gamma)$
generated by $A$ and   $\ttau_{c}$ where $1\le c< m+n$.

Then $\vphi_{w[n,m]}e(\gamma,\beta)$
belongs to $\tR$.

These generators satisfy the following commutation relations:
\eq&&\hs{4ex}\left\{\ba{l}
\tx_{a,b}\ttau_{c}
-\ttau_{c}\tx_{s_c(a),s_c(b)}\\[1ex]
\hs{7ex}=\hs{-6ex}\sum\limits_{\substack{
\nu_c=\nu_{c+1}\in\supp(\beta)\cap\supp(\gamma)}}
\hs{-5ex}\bl\delta(a=c+1)-\delta(a=c)-\delta(b=c+1)+\delta(b=c)\br e(\nu),\\[4ex]
\ttau_{a}^2
=\sum\limits_{\substack{\nu_a,\nu_{a+1}\in\supp(\beta)\cap\supp(\gamma)}}
Q_{\nu_a,\nu_{a+1}}(x_a,x_{a+1})e(\nu),\\[4ex]
\ttau_{a}\ttau_{b}-\ttau_{b}\ttau_{a}=0\quad\text{if $|a-b|>1$,}\\[2ex]
\ttau_{a+1}\ttau_{a}\ttau_{a+1}
-\ttau_{a}\ttau_{a+1}\ttau_{a}\\[1ex]
\hs{10ex}=
\hs{-5ex}\sum\limits_{\substack{
\nu_a,\nu_{a+1}\in\supp(\beta)\cap\supp(\gamma),\;\nu_a=\nu_{a+2}}}
\bQ_{\nu_a,\nu_{a+1}}(x_{a},x_{a+1},x_{a+2})e(\nu).
\ea\right.\label{eq:Arel}
\eneq
Indeed, the last equality follows from 
\begin{align*}
\ttau_{a+1}\ttau_{a}\ttau_{a+1}&=\sum\limits_{\nu}
\tau_{a+1}\tau_{a}\tau_{a+1}\;e(\nu)\quad\text{and}\\[.5ex]
\ttau_{a}\ttau_{a+1}\ttau_{a}&=\sum\limits_{\nu}
\tau_{a}\tau_{a+1}\tau_{a}\;e(\nu).
\end{align*}
Here the sums  in the both formulas range over
$\nu\in I^{\beta+\gamma}$ satisfying the conditions:
$\nu_a\in\supp(\gamma)$, $\nu_{a+1}\in\supp(\beta)\cap\supp(\gamma)$,
and $\nu_{a+2}\in \supp(\beta)$. 

Note that the error terms (i.e., the right hand sides of the equalities
in \eqref{eq:Arel})
belong to the algebra $A$ because we assume that 
$R(\beta)$ and $R(\gamma)$ are symmetric.
Hence we have
\eq&&\left\{\ba{l}
\tx_{a,b}\ttau_{c}
-\ttau_{c}\tx_{s_c(a),s_c(b)}\in A\\
\ttau_{a}^2\in A,\\
\ttau_{a}\ttau_{b}=\ttau_{b}\ttau_{a}\quad\text{if $|a-b|>1$,}\\
\ttau_{a+1}\ttau_{a}\ttau_{a+1}
-\ttau_{a}\ttau_{a+1}\ttau_{a}\in A.
\ea\right.\label{eq:sArel}
\eneq

Now for each element $w\in\sym_{m+n}$
let us choose  a  reduced expression $w=s_{a_1}\cdots s_{a_\ell}$.
We then set
$$\ttau_{w}=\ttau_{a_1}\cdots \ttau_{a_\ell}.$$
Then, similarly to a proof of the PBW decomposition \eqref{eq:PBW} (e.g., see \cite{KL09,R08}),
the commutation relations \eqref{eq:sArel} imply
$$\tR=\sum_{w\in\sym_{m+n}}\ttau_{w}A.$$
In particular, we obtain
$$\tR\subset \soplus_{\substack{w\in\sym_{n,m},\\w_1\in\sym_n,w_2\in\sym_m}}\tau_w \bl \tau_{w_1}\tens\tau_{w_2}\br
A.$$
It immediately implies \eqref{eq:-}, because we have,
for $1\le a<b\le m+n$, $\nu\in I^\gamma$, $\mu\in I^\beta$, $v\in e(\nu)N$ and $u\in e(\mu)M$,
\eqn
&&\tx_{a,b}\bl(v)_{z_2}\tens (u)_{z_1}\br\\
&&=\begin{cases}
\bl(x_a-x_b)v\br_{z_2}\tens (u)_{z_1}
&\hs{-31ex}\text{if $1\le a<b\le n$ and $\nu_a,\nu_b\in\supp(\beta)$,}\\[2ex]
(z_2-z_1)\bl(v)_{z_2}\tens (u)_{z_1}\br
+(x_av)_{z_2}\tens (u)_{z_1}-(v)_{z_2}\tens 
(x_{b-n}u)_{z_1}\\[.5ex]
&\hs{-31ex}\text{if $1\le a\le n<b\le m+n$ and}\\
&\hs{-15ex}\text{$\nu_a\in\supp(\beta),\mu_{b-n}\in\supp(\gamma)$,}\\[2ex]
(v)_{z_2}\tens \bl (x_{a-n}-x_{b-n})u\br_{z_1}
&\hs{-31ex}\text{if $n< a<b\le m+n$ and $\mu_{a-n},\mu_{b-n}\in\supp(\gamma)$,}
\\[1ex]
0&\hs{-31ex}\text{otherwise,}
\end{cases}
\eneqn
and $\bl \tau_{w_1}\tens\tau_{w_2}\br \bl(v)_{z_2}\tens (u)_{z_1}\br
=\bl\tau_{w_1}v\br_{z_2}\tens \bl\tau_{w_2}u\br_{z_1}$.

\end{document}